\documentclass[smallcondensed]{svjour3}     
\smartqed   
\usepackage{graphicx}
\usepackage{natbib}
\usepackage{amsmath}
\usepackage{amssymb}

\usepackage{epstopdf}

\newtheorem{thm}{Theorem}
\newtheorem{rmk}{Remark}
\newtheorem{lem}{Lemma}
\newtheorem{prop}{Proposition}
\newtheorem{defin}{Definition}
\newtheorem{expl}{Example}
\newtheorem{cor}{Corollary}
\journalname{Methodology and Computing in Applied Probability}
\begin{document}
\bibliographystyle{plainnat}

\title{Convergence in $L_p([0,T])$ of wavelet expansions of
$\varphi$-sub-Gaussian random processes}

\titlerunning{Wavelet expansions of
$\varphi$-sub-Gaussian processes}        

\author{Yuriy Kozachenko  \and  Andriy Olenko \and Olga~Polosmak }


\institute{Yu. Kozachenko \at
Department of Probability Theory,
Statistics and Actuarial Mathematics, Kyiv University, Kyiv,
Ukraine \\
              Tel.: +38-044-2590392\\
              Fax: +38-044-2590392\\
              \email{ykoz@ukr.net}          
                  \and
           A. Olenko \at
              Department of Mathematics and Statistics,
              La Trobe University,  Victoria, 3086, Australia\\
               Tel.: +61-3-94792609 \\
                            Fax: +61-3-94792466\\
                            \email{a.olenko@latrobe.edu.au}
\and
           O. Polosmak \at
Department of Economic Cybernetics, Kyiv University,
Kyiv, Ukraine\\
               Tel.: +38-044-2597043 \\
                            Fax: +38-044-2590392\\
                            \email{didenkoolga@yandex.ru}                           
}

\date{*Will appear in Methodology and Computing in Applied Probability.  The final publication is available at link.springer.com. DOI: 10.1007/s11009-013-9346-7.}

\maketitle

\begin{abstract}
The article presents new results on convergence in $L_p([0,T])$ of wavelet expansions of $\varphi$-sub-Gaussian random processes.  The  convergence rate of the expansions is obtained.  Specifications of the obtained results are discussed. 
\end{abstract}

\vspace{5cm}
\section{Introduction}
Multiresolution analysis of  deterministic signals by the use of wavelets has been extensively studied in recent years. However, in the context of stochastic processes, general wavelet approximations has not yet been fully investigated.  In the majority of cases developed deterministic methods and used error measures and metrics may not be appropriate to investigate wavelet representations of stochastic processes. It indicates the necessity of elaborating special stochastic techniques.

From a practical point of view, multi\-resolution analysis provides an efficient framework for the decomposition of random processes. 
Wavelet representations could be used to convert the problem of analyzing a continuous-time random process to that of analyzing a random sequence, which is much simpler. This approach is widely used in statistics to estimate a curve given observations of the curve plus some noise, in time series analysis for smoothing functional data, in simulation studies of various functionals defined on realizations of a random process, etc.

Recently, a considerable attention was given to wavelet orthonormal series representations of stochastic processes. Numerous results, applications, and references on convergence of wavelet expansions of random processes in various spaces  can be found in \citealt{att, bar, cam, cla, did, kur, koz2, koz3, koz5, zha}, just to mention a few. 

Figures~1 and 2\footnote{The figures have been generated by the R packages \textsc{geoR} and \textsc{wmtsa}. The Daubechies D8 wavelet basis and resolution levels 4 and 6 were used.} illustrate wavelet expansions of stochastic processes. A simulated realization of the process $\mathbf{X}(t)$ and its two wavelet reconstructions with different numbers of terms are plotted in Figure~1. Figure~2 displays boxplots of mean-square approximation errors for 500 simulated realizations of $\mathbf{X}(t)$ for each reconstruction.
Figure~2 suggests that empirical probabilities of large errors become smaller when the number of terms in the wavelet expansion increases.

Although the mentioned effect is well-known for deterministic functions, it has to be established theoretically for different stochastic processes and probability metrics. Numerical simulation results need to be confirmed by theoretical analysis. It is also important to obtain theoretical estimations of the rate of convergence for various stochastic wavelet expansions.

\begin{figure}[h]\label{fig1}  
\begin{center}
  \includegraphics[width=0.7\textwidth]{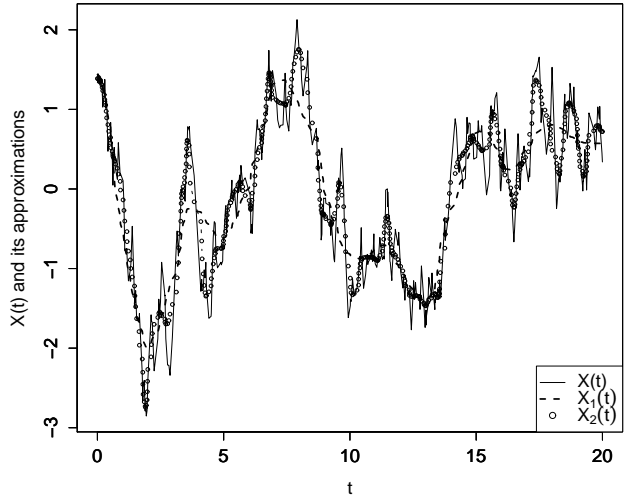}
\end{center}
\hspace{2cm}\parbox{10cm}{\caption{Stochastic process and its two wavelet reconstructions}}
\end{figure}

Our focus in this paper is on convergence in $L_p([0,T])$ of wavelet expansions of $\varphi$-sub-Gaussian random processes. The paper extends the recent results on uniform convergence  in the papers \citealt{koz2,koz3,koz5}  to new classes of stochastic processes and probability metrics.

The analysis and the approach presented in the paper contribute to the investigations of wavelet expansions of random processes in the former literature. The approach is of a special interest if $p > 2,$ as it
extends the available $L_2$ results.
In that sense, Theorems~\ref{conv_Lp_Sub_ner}-\ref{str_Sub_covLp_nonstation} are
of special importance. The results are obtained under simple assumptions which can be easily verified.  The paper deals with the most general class of wavelet expansions in comparison with particular cases considered by different authors, see, for example, \citealt{cam,kur}. 

The organization of the article is the following. In the second section we introduce the necessary background from the theory of $\varphi$-sub-Gaussian random variables and processes. Section 3 discusses wavelet expansions and approximations of stochastic processes. 
In \S 4 we present the main results on convergence in $L_p([0,T])$ of  wavelet expansions of $\varphi$-sub-Gaussian random processes. In this section we also obtain the rate of convergence of the wavelet expansions and discuss some specifications for which the assumptions in the theorems can be easily verified. 
\begin{figure}[t]\label{fig2}  
\begin{center}
  \includegraphics[width=0.7\textwidth]{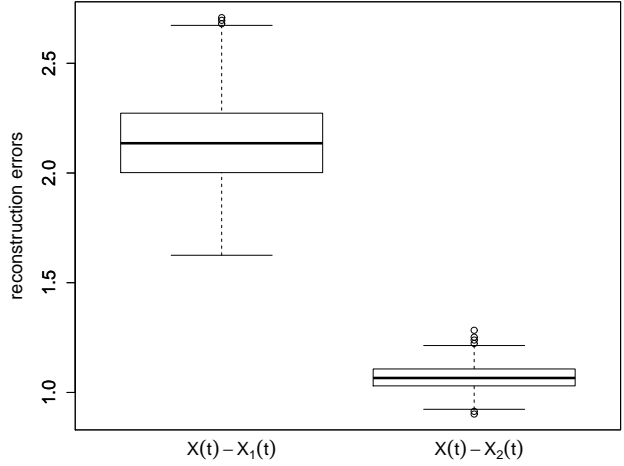}
\end{center}
\hspace{2cm}\parbox{10cm}{\caption{Boxplots of reconstruction errors for two approximations}}
\end{figure}
\section{$\varphi$-Sub-Gaussian random processes}
In this section, we review the definition of $\varphi$-sub-Gaussian random processes and  their relevant properties. 

The space of $\varphi$-sub-Gaussian random variables  was first introduced in the paper
\citealt{koz4}.  More information about the space of
$\varphi$-sub-Gaussian random variables and processes can be found in 
\citealt{bul,ant,koz0}.
\begin{defin} {\rm (\citealt{bul})} A continuous even convex function $\varphi(x),$ $x\in {\mathbb R},$ is called
an Orlicz N-function, if it is monotonically increasing for $x>0$,
$\varphi(0)=0,$ $ {{\varphi(x)}/{x}}\to 0,$ when $x\to 0,$ and $ {{\varphi(x)}/{x}}\to\infty,$  when $x\to\infty. $
\end{defin}

 \noindent{\bf Assumption Q} 
Let $\varphi(\cdot)$ be an Orlicz N-function and
$\lim\limits_{x\rightarrow 0}{\varphi(x)}/{x^2}=c>0.$

\begin{rmk} The constant $c$ can be equal to $+\infty.$ \end{rmk}

\begin{prop} {\rm(\citealt{bul})} Every N-function $\varphi(\cdot)$ can be represented as
$\varphi(u)=\int_{0}^{|u|}f(v)\,dv,$ where $f(\cdot)$ is
a monotonically nondecreasing, right-continuous function, such that
$f(0)=0$ and $f(x)\to\infty,$ when  $x\to\infty.$ The function 
$f(\cdot)$ is called a density of $\varphi(\cdot).$ \end{prop}

\begin{defin}{\rm (\citealt{bul})} Let $\varphi(x),x\in {\mathbb R},$ be an Orlicz N-function. The
function $\varphi^{*}(x),$ $x\in {\mathbb R},$ defined by the
formula $\varphi^{*}(x):=\sup_{y\in {\mathbb R}}(xy-\varphi(y))$ is called the
Young-Fenchel transform of $\varphi(\cdot).$ \end{defin}

Let $\{\Omega, \cal{B}, P\}$ be a standard probability space and $L_p(\Omega)$ denote a space of random variables having finite $p$-th absolute moments.

\begin{defin} {\rm(\citealt{ant})} Let $\varphi(\cdot)$ be an Orlicz
N-function satisfying the assumption $Q.$ A zero mean random
variable $\xi$ belongs to the space $Sub_\varphi(\Omega)$ (the space
of $\varphi$-sub-Gaussian random variables), if there exists a
constant $r_{\xi}\geq 0$ such that the inequality $\mathbf
E\exp\left(\lambda\xi\right)\le
\exp\left(\varphi(r_{\xi}\lambda)\right)$ holds for all $\lambda\in
{\mathbb R}.$ \end{defin}

\begin{prop} {\rm (\citealt{koz4})} 
The space $Sub_\varphi(\Omega)$ is a Banach space with respect to
the norm $$ \tau_\varphi(\xi):=\inf\{a\ge0: \mathbf
E\exp\left(\lambda\xi\right)\le \exp(\varphi(a\lambda)),\lambda\in
{\mathbb R}\}.$$\end{prop}

\begin{rmk} A Gaussian centered random variable $\xi$ belongs to the space $Sub_\varphi(\Omega),$
where $\varphi(x)={x^2}/2$ and $\tau^2_\varphi(\xi)=\mathbf
E\xi^2.$
\end{rmk}

\begin{defin} {\rm (\citealt{bul})} Let $\mathbf{T}$ be a parametric space. A random process $\mathbf
X(t),$ $t \in \mathbf{T},$ belongs to the space $Sub_{\varphi}(\Omega)$ if $\mathbf X(t)\in Sub_{\varphi}(\Omega)$ for
all $t\in \mathbf{T}$  and $\sup_{t\in
\mathbf{T}}\tau_\varphi(\mathbf X(t))<\infty.$
\end{defin}

Let us denote $\tau_{\varphi}(t):=\tau_{\varphi}\left(\mathbf X(t)\right).$

\begin{defin} {\rm(\citealt{koz1})} A family $\Xi$ of random
variables $\xi\in Sub_{\varphi}(\Omega)$ is called strictly
$Sub_{\varphi}(\Omega)$ if there exists a constant $C_{\Xi}>0$
such that for any finite set $I,$ $\xi_{i}\in \Xi,$ $i\in I,$
and for arbitrary $\lambda_{i}\in {\mathbb R},$ $i\in I:$ 
\[\tau_{\varphi}\left(\sum\limits_{i \in
I}\lambda_{i}\xi_{i}\right)\leq
 C_{\Xi}\left(\mathbf E\left(\sum\limits_{i \in I}\lambda_{i}\xi_{i}\right)^2\right)^{1/2}.\]
$C_{\Xi}$ is called a determinative constant. \end{defin}

\begin{defin}\label{def_strictly Sub} {\rm(\citealt{koz1})} $\varphi$-sub-Gaussian random
process $\mathbf X(t),$ $t\in \mathbf{T},$ is called strictly
$Sub_{\varphi}(\Omega)$ if the family of random variables $\{\mathbf
X(t),t\in \mathbf{T}\}$ is strictly $Sub_{\varphi}(\Omega).$ The determinative
constant of this family is called a determinative constant of the
process and denoted by $C_{X}$.
\end{defin}

\begin{rmk} Gaussian centered random process $\mathbf X(t),$ $t\in \mathbf{T},$ is a strictly
$Sub_{\varphi}(\Omega)$ process, where $\varphi(x)={x^2}/2$ and the
determinative constant $C_{X}=1.$
\end{rmk}

\begin{thm}\label{ner_XSub_exp}{\rm(\citealt{koz0})}  Let $\{\mathbf{T},\Lambda,\mu\}$ be a measurable space and $\mathbf X(t),$ $t\in \mathbf{T},$
be a random process from the space $Sub_{\varphi}(\Omega).$

If $\int_\mathbf{T}
\left(\tau_{\varphi}(t)\right)^p\,d\mu(t)<\infty$ for some $p\ge 1,$ then the integral $\int_\mathbf{T}
\left|\mathbf X(t)\right|^p\,d\mu(t)$ exists with probability {\rm 1} and the following
inequality holds

\[\mathbf P\left\{\int_\mathbf{T}
\left|\mathbf
X(t)\right|^p\,d\mu(t)>\varepsilon\right\}\le 2\exp\left\{-\varphi^*\left(\left({\varepsilon}/c\right)^{1/p}\right)\right\},\]
for each non-negative
\begin{equation}\label{epsilon_ner}\varepsilon>c\cdot\left(f\left(p(c/\varepsilon)^{1/p}\right)\right)^p,\end{equation}
where $c:=\int_\mathbf{T} \left(\tau_{\varphi}(t)\right)^p\,d\mu(t)$
and $f(\cdot)$ is a density of $\varphi(\cdot).$ 

\end{thm}

\begin{expl} If $\varphi(x)={|x|^{\alpha}}/\alpha,$ $
1<\alpha\le2,$ then $f(x)=x^{\alpha-1}$ and
$\varphi^*(x)={|x|^{\beta}}/\beta,$ where $\beta>1$ and $1/{\alpha}+1/{\beta}=1.$ Hence,
$c\cdot\left(f\left(p(c/\varepsilon)^{1/p}\right)\right)^p=c^{\alpha}p^{\alpha-1}\varepsilon^{1-\alpha}.$
Therefore, inequality {\rm(\ref{epsilon_ner})} holds if
$\varepsilon>c\cdot p^{\frac{\alpha-1}{\alpha}p}.$
\end{expl}

\begin{rmk} If $\mathbf X(t),$ $t\in \mathbf{T},$ is a Gaussian centered random
process and $\sigma(t):=\left(\mathbf E\left(\mathbf
X(t)\right)^2\right)^{1/2},$ then inequality {\rm(\ref{epsilon_ner})} holds
true for $\varepsilon>\tilde c\cdot p^{\frac{p}{2}},$ where $\tilde
c:=\int_\mathbf{T} \left(\sigma(t)\right)^p\,d\mu(t).$ For such
$\varepsilon$ the following estimate is valid
\[\label{intX(t)gauss_ner}\mathbf P\left\{\int_\mathbf{T}
\left|\mathbf
X(t)\right|^p\,d\mu(t)>\varepsilon\right\}\le2\exp\left\{-\frac12\left(\frac{\varepsilon}{\tilde
c}\right)^{2/p}\right\}.\]
\end{rmk}

\begin{rmk} If $\mathbf X(t),$ $t\in \mathbf{T},$ is a strictly
$Sub_{\varphi}(\Omega)$ random process, then
\[c\le C^p_{X}\int_\mathbf{T}
\left(\mathbf E\left(\mathbf X(t)\right)^2\right)^{\frac
p2}\,d\mu(t).\]
\end{rmk}

\section{Wavelet representation of random processes}
In this section we introduce wavelet representations and approximations of non-random functions and stochastic processes.  

Let $\phi(x),$ $x\in\mathbb R,$ be a function from the space
$L_2(\mathbb R)$ such that $\widehat{\phi}(0)\ne 0$ and
$\widehat{\phi}(y)$ is continuous at $0,$ where
$\widehat{\phi}(y)=\int_{\mathbb R}e^{-iyx}{\phi(x)}\,dx$ is the
Fourier transform of $\phi.$

 Suppose that the following assumption
holds true: $\sum_{k\in\mathbb Z} |\widehat{\phi}(y+2{\pi}k)|^2=1\
{\rm (a.e.)} $

There exists a function $m_0(x)\in L_2([0,2\pi])$, such that
$m_0(x)$ has the period $2\pi$ and
$$\widehat{\phi}(y)=m_0\left(y/2\right)\widehat{\phi}\left(y/2\right)\ {\rm (a.e.)}
$$
 In this case the  function $\phi(\cdot)$ is called
the $f$-wavelet.

Let $\psi(\cdot)$ be the inverse Fourier transform of the function
$$\widehat{\psi}(y):=\overline{m_0\left(\frac
y2+\pi\right)}\cdot\exp\left\{-i\frac
y2\right\}\cdot\widehat{\phi}\left(\frac y2\right).$$ Then the
function $\psi(x)=\frac1{2\pi}\int_{\mathbb
R}e^{iyx}{\widehat{\psi}(y)}\,dy$ is called the $m$-wavelet.

Let
\[\phi_{jk}(x):=2^{j/2}\phi(2^jx-k),\quad
\psi_{jk}(x):=2^{j/2}\psi(2^jx-k),\quad j,k \in\mathbb Z\,.
\]

It is known that the family of functions $\{\phi_{0k},
\psi_{jk}:j\in \mathbb N_0,\,k\in \mathbb Z\}$ is an orthonormal basis in
$L_2(\mathbb R)$ (see, for example, \citealt{dau}).

An arbitrary function $f(x)\in L_2(\mathbb R)$ can be represented in
the form
\begin{equation}\label{2.5}f(x)=\sum_{k\in\mathbb Z}\alpha_{0k}\phi_{0k}(x)+\sum_{j=0}^{\infty}\sum_{k\in\mathbb Z}\beta_{jk}\psi_{jk}(x)\,,
\end{equation}
$$\alpha_{0k}:=\int_{\mathbb R}f(x)\overline{\phi_{0k}(x)}\,dx,\quad \beta_{jk}:=\int_{\mathbb R}f(x)\overline{\psi_{jk}(x)}\,dx.$$
The representation (\ref{2.5}) is called a wavelet representation.

The series~(\ref{2.5}) converges in $L_2(\mathbb R)$ i.e.
$\sum_{k\in\mathbb
Z}|\alpha_{0k}|^2+\sum_{j=0}^{\infty}\sum_{k\in\mathbb
Z}|\beta_{jk}|^2<\infty\,.$

The integrals $\alpha_{0k}$ and $\beta_{jk}$ may also exist for
functions from $L_p(\mathbb R)$ and other function spaces. Therefore
it is possible to obtain the representation~(\ref{2.5}) for function
classes which are wider than $L_2(\mathbb R)$ (see, for example,
\citealt{jaf,tri}).

\noindent {\bf Assumption S} (\citealt{har})  For the
function $\phi(\cdot)$ there exists a decreasing function $\Phi(x),$ $x\ge
0,$ such that $\Phi(0)<\infty,$ $|\phi(x)|\le \Phi(|x|)$ (a.e.), and the integral
$\displaystyle\int_{\mathbb R}\Phi(|x|)\,dx$ is finite.

 Let
$\mathbf{X}(t),$ $t\in\mathbb R,$ be a  random process such that
$\mathbf E\mathbf{X}(t)=0\,$ for all $t\in\mathbb R$.
If sample trajectories of this process are  in the space
$L_2(\mathbb R)$ with probability one, then it is possible to obtain
the representation (wavelet representation)
\begin{equation}\label{X(t)}\mathbf{X}(t)=\sum_{k\in \mathbb Z}\xi_{0k}\phi_{0k}(t)+\sum_{j=0}^{\infty}\sum_{k\in \mathbb Z}
\eta_{jk}\psi_{jk}(t)\,, \end{equation} with mean-square integrals
\begin{equation}\label{intksieta}\xi_{0k}:=\int_{\mathbb R}\mathbf{X}(t)\overline{\phi_{0k}(t)}\,dt,\quad
\eta_{jk}:=\int_{\mathbb
R}\mathbf{X}(t)\overline{\psi_{jk}(t)}\,dt\,.\end{equation}

However, the majority of random processes does not possess the required
property. For example, sample paths of stationary processes are not
in the space $L_2(\mathbb R)$ (a.s.). However, in many cases it is
possible to construct a representation of type~(\ref{X(t)}) for
$\mathbf{X}(t).$

Consider the approximants of $\mathbf{X}(t)$ defined by
\[\label{x_nk}
\mathbf{X}_{n,\mathbf{k}_n}(t):=\sum_{|k|\le
k_0'}\xi_{0k}\phi_{0k}(t)+\sum_{j=0}^{n-1}\sum_{|k|\le
k_j}\eta_{jk}\psi_{jk}(t)\,,
\]
where $\mathbf{k}_n:=(k_0'(n),k_0(n),...,k_{n-1}(n)).$

For the sake of simplicity, we will omit the index $(n)$ in $(k_0'(n),k_0(n),...,k_{n-1}(n))$ and will use the notation $\mathbf{k}_n\to \infty$ to denote that  $n\to\infty,$ $k_0'\to\infty,$ and $k_j\to\infty$ for all
$j\in \mathbb N_0\,.$ 

Theorem~\ref{8103} below guarantees the mean-square convergence of
$\mathbf{X}_{n,\mathbf{k}_n}(t)$  to $\mathbf{X}(t).$

\begin{thm}\label{8103} {\rm (\citealt{koz5})} Let $\mathbf{X}(t),$ $t\in\mathbb R,$ be a random process
such that $\mathbf E\mathbf{X}(t)=0,$ $\mathbf
E|\mathbf{X}(t)|^2<\infty$ for all $t\in\mathbb R,$ and its
covariance function  $R(t,s)$  is continuous. Let the $f$-wavelet
$\phi$ and the $m$-wavelet $\psi$ be continuous functions and  the
assumption {\rm S}  hold true for both  $\phi$ and  $\psi.$  Suppose that there exist a function
$A:(0,\infty)\to (0,\infty)$ and $x_0\in\mathbb R$ such that  $
c(ax)\le c(x)\cdot A(a),$ for all $x\ge x_0.$

Assume that there exists an even non-decreasing  on $[0,\infty)$  function $c(x),$ $x\in\mathbb R,$
 with $c(0)>1$ such that $\int_{\mathbb R}c(x)\Phi(|x|)\,dx<\infty\quad {\rm and}
\quad |R(t,t)|^{1/2}\le c(t)$ for all $t\in\mathbb R.$ Then
\begin{enumerate}
  \item[\rm 1.]  $\mathbf{X}_{n,\mathbf{k}_n}(t)\in L_2(\Omega)\,;$
  \item[\rm 2.]  $\mathbf{X}_{n,\mathbf{k}_n}(t)\to \mathbf{X}(t)$ in mean square when  $\mathbf{k}_n\to\infty.$
\end{enumerate}
\end{thm}

\begin{rmk} It was shown  that the integrals in {\rm(\ref{intksieta})} exist if $\mathbf{X}(t)$ satisfies the assumptions of 
Theorem~{\rm\ref{8103},} see {\rm \citealt{koz5}}.  Though there are other sufficient conditions in the literature which guarantee the existence of the integrals in {\rm(\ref{intksieta})} and the convergence in the mean-square sense in Theorem~{\rm\ref{8103},} the assumptions in terms of the functions $c(\cdot)$ and $\Phi(\cdot)$  can be easily verified in many practical examples.
\end{rmk}

\section{Convergence of wavelet expansions of
$\varphi$-sub-Gaussian random processes}

In this section we present the main results on convergence in $L_p(\mathbf{T}),$ $\mathbf{T}=[0,T],$ $T>0,$ of the wavelet expansions of $\varphi$-sub-Gaussian random processes. The rate of convergence in the space $L_p([0,T])$ is obtained. We also present some specifications of the general results for which the assumptions can be easily verified. 

In \citealt{ant}, it was shown that there always exists a constant $c_{\varphi}>0$
such that $\mathbf E|\xi|^2\le c_{\varphi}\tau^2_{\varphi}(\xi).$ Therefore,  random processes from the space
$Sub_{\varphi}(\Omega)$ belong to the space
$L_{2}(\Omega).$ The following theorem is a corollary of  Theorem~\ref{8103} and the estimate $\left(R(t,t)\right)^{1/2}\le
\sqrt{c_{\varphi}}\tau_{\varphi}(t).$

\begin{thm}\label{meansq_subfi} Let $\mathbf{X}(t),$ $t\in\mathbb R,$ be a random process
such that $\mathbf{X}(\cdot)\in Sub_{\varphi}(\Omega)$ and $\mathbf{X}(t)$
is continuous in the norm $\tau_{\varphi}(\cdot).$ Let the functions 
 $\phi(\cdot),$  $\psi(\cdot),$  $\Phi(\cdot),$ $A(\cdot),$ and $c(\cdot)$ satisfy the assumptions of Theorem~{\rm\ref{8103}}. 
 
If $\int_{\mathbb R}c(x)\Phi(|x|)\,dx<\infty\quad {\rm and}
\quad \tau_{\varphi}(t)\le c(t)$ for all $t\in\mathbb R,$ then
\begin{enumerate}
  \item[\rm 1.]  $\mathbf{X}_{n,\mathbf{k}_n}(t)\in L_2(\Omega)\,;$
  \item[\rm 2.]  $\mathbf{X}_{n,\mathbf{k}_n}(t)\to \mathbf{X}(t)$ in mean square when  $\mathbf{k}_n\to\infty.$
\end{enumerate}
\end{thm}

\begin{defin} A random
process $\mathbf X(t),$ $t\in \mathbb{R},$ belongs to the space $L_p([0,T])$ if

\[\mathbf P \left\{\left(\int_0^T\left|\mathbf
X(t)\right|^p\,dt\right)^{1/p}<+\infty\right\}=1.\]
\end{defin}

\begin{thm}\label{conv_Lp_Sub_ner} Let $\mathbf{X}(t),$ $t\in\mathbb R,$ be a random process
such that $\mathbf{X}(\cdot)\in Sub_{\varphi}(\Omega)$ and $\mathbf{X}(t)$
is measurable and continuous in the norm $\tau_{\varphi}(\cdot).$
Let the assumptions of Theorem~{\rm\ref{meansq_subfi}} are satisfied.
Suppose that 
\begin{equation}\label{inttau}\int_0^T\left(\tau_{\varphi}(t)\right)^p\,dt<\infty\end{equation} for
some $T>0$ and $p\ge 1.$ 

If 
\begin{equation}\label{inttau_Xn-Xm}\int_0^T\left(\tau_{\varphi}\left(\mathbf X_{n,\mathbf k_n}(t)
-\mathbf X_{m,\mathbf k_m}(t)\right)\right)^p\,dt\to 0,\end{equation}
when $\mathbf k_n,\mathbf k_m\to\infty,$  then $\mathbf X_{n,\mathbf k_n}(t)\to\mathbf X(t)$ in the space $L_p([0,T])$ in
probability, when $\mathbf k_n\to\infty,$ i.e.  for all $\varepsilon>0$
\begin{eqnarray*}\label{P_int_Xn-X}\mathbf P \left\{\left(\int_0^T\left|\mathbf
X_{n,\mathbf k_n}(t) -\mathbf
X(t)\right|^p\,dt\right)^{1/p}>\varepsilon\right\}\to 0,\end{eqnarray*}
when $\mathbf k_n\to\infty$. 

Furthermore, the following estimate
holds
\begin{equation}\label{P_int_Xn-Xm<exp}\mathbf P \left\{\int_0^T\left|\mathbf
X_{n,\mathbf k_n}(t) -\mathbf
X(t)\right|^p\,dt>\varepsilon\right\}\le2\exp\left\{-{\varphi}^*\left(\left({\varepsilon}/{c_n^{\infty}}
\right)^{1/p}\right)\right\},\end{equation} where
$c_n^{\infty}:=\int_0^T\left(\tau_{\varphi}\left(\mathbf
X_{n,\mathbf k_n}(t) -\mathbf X(t)\right)\right)^p\,dt$ and
$\varepsilon>c_n^{\infty}\left(f\left(p\,(c_n^{\infty}/\varepsilon)^{1/p}
\right)\right)^p.$
\end{thm}
\noindent{\it Proof.} First, we will show that the random variables $\xi_{0k}$
and $\eta_{jk}$ are from the space $Sub_{\varphi}(\Omega).$ We prove it only for $\eta_{jk}.$ The case of $\xi_{0k}$ can be dealt with similarly.

We will need the following generalization of the Minkowski inequality.
\begin{lem}\label{gmink}
If $\mathbf{X}(t),$ $t\in\mathbb R,$ is a measurable $\varphi$-sub-Gaussian random processes, then
\[\tau_{\varphi}\left(\int_{\mathbb R}\mathbf{X}(t)\,dt\right)\le\int_{\mathbb R}\tau_{\varphi}(t)\,dt.\]
\end{lem}
\noindent{\it Proof.} Without loss of generality we may assume that $\tau_{\varphi}(t)>0$ for all $t\in\mathbb R$  and $\int_{\mathbb R}\tau_{\varphi}(t)\,dt<\infty.$ 

Let us denote
\begin{eqnarray*}
 I: &=&  \mathbf E \left( \exp \left\{  \lambda
\int_{\mathbb R} \mathbf{X}(t)\,dt \cdot \left(
\int_{\mathbb R} \tau_{\varphi}(s)\,ds \right)^{-1} 
\right\}\right)\\
&=& \mathbf E \left(\exp \left\{  \lambda
\int_{\mathbb R} \frac{\mathbf{X}(t)}{\tau_{\varphi}(t)} \cdot
\frac{\tau_{\varphi}(t)}{\int_{\mathbb R} \tau_{\varphi}(s)\,ds}\,dt
 \right\}\right),\ 
\lambda \in \mathbb R.
\end{eqnarray*}

Notice that $\exp(\cdot)$ is a convex function and
\[
  \int_{\mathbb R} \tau_{\varphi}(t) \cdot \left(
\int_{\mathbb R} \tau_{\varphi}(s)\,ds \right)^{-1} \,dt =1.
\]

Hence, by Jensen's inequality, we obtain
\[
\exp \left\{ \lambda
\int_{\mathbb R} \frac{\mathbf{X}(t)}{\tau_{\varphi}(t)} \cdot
\frac{\tau_{\varphi}(t)}{\int_{\mathbb R} \tau_{\varphi}(s)\,ds}\,dt
 \right\} \leq \int_{\mathbb R}  \exp \left\{ \lambda
\frac{\mathbf{X}(t)}{\tau_{\varphi}(t)}\right\} \cdot
\frac{\tau_{\varphi}(t)}{\int_{\mathbb R} \tau_{\varphi}(s)\,ds}\,dt.
\] 

Therefore,
\begin{equation}\label{I}
I \leq \int_{\mathbb R} \mathbf E \left(\exp \left\{ \lambda
\frac{\mathbf{X}(t)}{\tau_{\varphi}(t)}\right\}\right) \cdot
\frac{\tau_{\varphi}(t)}{\int_{\mathbb R} \tau_{\varphi}(s)\,ds}\,dt\,.
\end{equation}

By the definition of $\tau_{\varphi}(\cdot)$ we get
\begin{equation}\label{est}
 \mathbf E\left( \exp \left\{ \lambda
\frac{\mathbf{X}(t)}{\tau_{\varphi}(t)}\right\}\right) \leq \exp \{ \varphi (\lambda) \}\,.
\end{equation}

By inequalities (\ref{I}) and (\ref{est}) we obtain $I\leq \exp \{ \varphi (\lambda) \}.$ Finally,  the statement of the lemma follows from the definition of $\tau_{\varphi}(\cdot).$
\hfill   $\Box$

By Lemma~\ref{gmink} we get
\begin{eqnarray*}\tau_{\varphi}(\eta_{jk})&\le&\int_{\mathbb R}\tau_{\varphi}(t)\left|\psi_{jk}(t)\right|\,dt\le
\int_{\mathbb R}\tau_{\varphi}(t)2^{j/2}\left|\psi(2^j
t-k)\right|\,dt\\ 
&\le& 2^{j/2} \int_{\mathbb R}c(t)|\psi(2^j t-k)|\,dt=2^{-j/2}\int_{\mathbb R}c\left(\frac{|u+k|}{2^j}\right)\psi(u)\,du\\
&\le& 2^{-j/2} \int_{\mathbb R}c\left(\frac{|u+k|}{2^j}\right)\Phi(|u|)\,du<+\infty,
\end{eqnarray*}
since $c(|u+k|)\le c(2u)\le c(u)A(2)$ for sufficiently large $u$.

Thus, the processes $\mathbf X_{n,\mathbf k_n}(t)$ belongs to the
space $Sub_{\varphi}(\Omega).$

To prove the theorem, it is enough to show that for all
$\varepsilon>0$
\begin{eqnarray}\label{P_int_Xn-Xm1}\mathbf P \left\{\int_0^T\left|\mathbf
X_{n,\mathbf k_n}(t) -\mathbf X_{m,\mathbf
k_m}(t)\right|^p\,dt>\varepsilon\right\}\to 0,\end{eqnarray} when
$\mathbf k_n$ and $\mathbf k_m\to+\infty.$

Indeed, by Theorem~\ref{ner_XSub_exp} and (\ref{inttau}), $\mathbf X(t)\in
L_p([0,T])$ with probability 1.

If the condition (\ref{P_int_Xn-Xm1}) is satisfied then
there exists a process $\mathbf Y(t)\in L_p([0,T])$ such that
\begin{eqnarray*}\mathbf P \left\{\int_0^T\left|\mathbf
X_{n,\mathbf k_n}(t) -\mathbf
Y(t)\right|^p\,dt>\varepsilon\right\}\to 0,\end{eqnarray*} 
when
$\mathbf k_n\to+\infty,$
and by Theorem \ref{meansq_subfi} we get $\mathbf Y(t)=\mathbf X(t)$ with probability 1. 

To prove (\ref{P_int_Xn-Xm1}) we denote
\begin{eqnarray*}c^m_n:=\int_0^T\left(\tau_{\varphi}\left(\mathbf X_{n,\mathbf k_n}(t)
-\mathbf X_{m,\mathbf k_m}(t)\right)\right)^p\,dt.\end{eqnarray*} 
By Theorem \ref{ner_XSub_exp} the following inequality holds
\begin{equation}\label{intXn-Xm_ner_cnm}\mathbf P\left\{\int _T
\left|\mathbf X_{n,\mathbf k_n}(t) -\mathbf X_{m,\mathbf
k_m}(t)\right|^p\,dt>\varepsilon\right\}\le 2\exp\left\{-\varphi^*\left(\left({\varepsilon}/{c^m_n}\right)^{1/p}\right)\right\}\end{equation}
for
$\varepsilon>c^m_n\cdot\left(f\left(p\,(c^m_n/\varepsilon)^{1/p}\right)\right)^p.$

The condition (\ref{inttau_Xn-Xm}) implies  that
$c^m_n\to 0$ as $n,m\to\infty.$ The function $f(x)$ is
right-continuous and $f(0)=0.$ Hence
$f\left(p\,(c^m_n/\varepsilon)^{1/p}\right)\to 0, \ \
n,m\to\infty.$ Therefore, for arbitrary $\varepsilon>0$ and  sufficiently large
$n,m$  inequality (\ref{intXn-Xm_ner_cnm}) holds true. By definitions~1 and 2 the right-hand side of the inequality vanishes when $c^m_n\to 0,$ which implies  (\ref{P_int_Xn-Xm1}).

Notice, that $\mathbf X_{n,\mathbf k_n}(t) -\mathbf X_{m,\mathbf
k_m}(t)\to \mathbf X_{n,\mathbf k_n}(t) -\mathbf X(t)$ in $L_p([0,T])$ with probability 1 and $c^m_n\to c^\infty_n,$ when $\mathbf k_m\to+\infty.$
Hence, we obtain inequality (\ref{P_int_Xn-Xm<exp}) if we let $\mathbf k_m$ in (\ref{intXn-Xm_ner_cnm}) approach infinity.  \hfill   $\Box$

\begin{cor}\label{cor1_str_Sub}
Let  $\mathbf{X}(t),$ $t\in\mathbb R,$ be a strictly
$Sub_{\varphi}(\Omega)$ random process with a determinative constant $C_{X}.$ Then the
 condition {\rm(\ref{inttau})} holds true if  the integral
$\int_0^T\left(\mathbf E(\mathbf X(t))^2\right)^{p/2}\,dt$ is convergent.

If
\begin{equation}\label{int_Xn-Xm^2}\int_0^T\left(\mathbf E\left(\mathbf X_{n,\mathbf k_n}(t)
-\mathbf X_{m,\mathbf
k_m}(t)\right)^2\right)^{p/2}\,dt\to 0,\end{equation}
then the condition {\rm(\ref{inttau_Xn-Xm})} holds when $\mathbf k_n,\mathbf k_m\to\infty.$
The inequality {\rm(\ref{P_int_Xn-Xm<exp})} is valid for 
\begin{equation}\label{c_n^infty}
c_n^{\infty}=C_{X}^p\int_0^T\left(\mathbf E\left(\mathbf
X_{n,\mathbf k_n}(t) -\mathbf
X(t)\right)^2\right)^{p/2}\,dt.\end{equation}
\end{cor}
\noindent{\it Proof.} The statement of  Corollary \ref{cor1_str_Sub}  follows from
 Theorem \ref{conv_Lp_Sub_ner}, Definition \ref{def_strictly Sub},
and the fact that a linear closure of a family of strictly
$Sub_{\varphi}(\Omega)$ random variables is a family of strictly
$Sub_{\varphi}(\Omega),$ see \citealt{koz1}. \hfill   $\Box$

\begin{cor}\label{cor2_str_Sub}
The statement of Corollary~{\rm\ref{cor1_str_Sub}} is valid if instead of  {\rm(\ref{int_Xn-Xm^2})} the  series 
\[\sum_{k\in \mathbb Z}\left(\mathbf E\xi_{0k}^2\right)^{1/2}\left|\phi_{0k}(t)\right|
+\sum_{j=0}^{\infty}\sum_{k\in \mathbb Z}\left(\mathbf
E\eta_{jk}^2\right)^{1/2}\left|\psi_{jk}(t)\right| \]
converges uniformly on the interval $[0,T].$
Then the inequality {\rm(\ref{P_int_Xn-Xm<exp})} holds for
\begin{eqnarray}\label{c_n^infty2}
c_n^{\infty}=C_{X}^p\ T\, \cdot&&\hspace{-2mm}\sup\limits_{0\le t \le
T}\left(\sum_{|k|\ge
 k'_0}\left(\mathbf
E\xi_{0k}^2\right)^{1/2}\left|\phi_{0k}(t)\right|+\sum_{j=0}^{n-1}\sum_{|k|\ge
 k_j}\left(\mathbf
E\eta_{jk}^2\right)^{1/2}\left|\psi_{jk}(t)\right|\right.\nonumber\\
& &+\left.\sum_{j=n}^{\infty}\sum_{k\in
\mathbb Z}\left(\mathbf
E\eta_{jk}^2\right)^{1/2}\left|\psi_{jk}(t)\right|\right)^p.\nonumber\end{eqnarray}
\end{cor}
\noindent{\it Proof.}  The statement of Corollary \ref{cor2_str_Sub}  follows from (\ref{c_n^infty}) and
the estimate
$$\left(\mathbf E\left(\mathbf X_{n,\mathbf k_n}(t)
-\mathbf X(t)\right)^2\right)^{1/2}=\left(\mathbf
E\left(\sum_{|k|\ge
 k'_0}\xi_{0k}\phi_{0k}(t)+\sum_{j=0}^{n-1}\sum_{|k|\ge
  k_j}\eta_{jk}\psi_{jk}(t)\right.\right.$$
  $$\left.\left.+\sum_{j=n}^{\infty}\sum_{k\in
  \mathbb Z}\eta_{jk}\psi_{jk}(t)\right)^2\right)^{1/2}\le\sum_{|k|\ge
 k'_0}\left(\mathbf
E\xi_{0k}^2\right)^{1/2}\left|\phi_{0k}(t)\right|+\sum_{j=0}^{n-1}\sum_{|k|\ge
 k_j}\left(\mathbf
E\eta_{jk}^2\right)^{1/2}$$
  $$\times\left|\psi_{jk}(t)\right|+\sum_{j=n}^{\infty}\sum_{k\in
\mathbb Z}\left(\mathbf
E\eta_{jk}^2\right)^{1/2}\left|\psi_{jk}(t)\right|.$$
 \hfill   $\Box$
\begin{lem}\label{lem_ner_phi} If $\delta(\cdot)$ is an m-wavelet or an f-wavelet satisfying the assumption $S,$ then
\begin{equation}\label{sum|delta|}
\sup\limits_{x\in \mathbb R}\,\sum\limits_{k\in\mathbb Z}|\delta(x-k)|\le
3\Phi(0)+4\int_{1/2}^{\infty}\Phi(t)\,dt=:C_{\delta}
\end{equation}
and  
\begin{equation}\label{supsum|delta|}\sup\limits_{|x|\le T}\sum\limits_{|k|\ge k_1}|\delta(x-k)|\le
\int_{k_1-T-1}^{\infty}\Phi(t)\,dt+\int_{k_1-1}^{\infty}\Phi(t)\,dt=:C_{\delta}(T,k_1)
\end{equation}
for $k_1\ge T+1.$
\end{lem}
\noindent{\it Proof.} Inequality (\ref{sum|delta|}) is a simple
modification of an inequality from \citealt{har}. Therefore we
only prove (\ref{supsum|delta|}).

By the assumption $S,$
$$\sum\limits_{|k|\ge k_1}|\delta(x-k)|\le \sum\limits_{k\ge k_1}\left(\Phi(|x+k|)
+\Phi(|x-k|)\right)=:z_{k_1}(x),$$ where $z_{k_1}(x)$ is an even
function. As $\Phi(\cdot)$ is a decreasing function on $[0,+\infty)$ and $k_1\ge T+1,$ we obtain
\begin{eqnarray*}\sup\limits_{|x|\le T}\sum\limits_{|k|\ge k_1}|\delta(x-k)|\le\sup\limits_{0\le x\le T}z_{k_1}(x)\le\sum\limits_{k\ge
k_1}\left(\Phi(|k-T|) +\Phi(|k|)\right)\\ \le\sum\limits_{k\ge k_1}
\left(\int_{k-1}^{k}\Phi(t-T)\,dt+\int_{k-1}^{k}\Phi(t)\,dt\right)\le
\int_{k_1-T-1}^{\infty}\Phi(t)\,dt+\int_{k_1-1}^{\infty}\Phi(t)\,dt.
\end{eqnarray*}
\hfill   $\Box$
\begin{thm}
Let $\mathbf{X}(t),$ $t\in\mathbb R,$ be a strictly
$Sub_{\varphi}(\Omega)$ random process for which the assumptions of
 Theorem~{\rm\ref{8103}} hold true and
$\int_0^T\left(\mathbf E(\mathbf X(t))^2\right)^{p/2}\,dt<+\infty.$

If 
\begin{equation}\label{supsumksieta^2}\sup_{k\in \mathbb Z}\left(\mathbf
E|\xi_{0k}|^2\right)^{1/2}C_{\phi} +\sum_{j=0}^{\infty}\sup_{k\in
\mathbb Z}\left(\mathbf
E|\eta_{jk}|^2\right)^{1/2}\cdot 2^{j/2}C_{\psi}<\infty,\end{equation}
 then $\mathbf X_{n,\mathbf k_n}(t) \to\mathbf X(t)$ in the space $L_p(0,T)$ in probability, when
$\mathbf k_n\to\infty.$

Furthermore, for $\varepsilon>c_n^{\infty}\left(f\left(p\,(c_n^{\infty}/\varepsilon)^{1/p}
\right)\right)^p$  inequality
{\rm (\ref{P_int_Xn-Xm<exp})} is valid when
\begin{eqnarray}
c_n^{\infty}=C_{X}^p\ T&&\hspace{-1mm}\left(\sup_{k\in \mathbb Z}\left(\mathbf
E|\xi_{0k}|^2\right)^{1/2}C_{\phi}(T,k_0')+\sum_{j=0}^{J-1}\sup_{k\in \mathbb Z}\left(\mathbf
E|\eta_{jk}|^2\right)^{1/2}\right.\nonumber\\ 
&&\times\left.2^{j/2}C_{\psi}(T,k_j)+\sum_{j=J}^{\infty}\sup_{k\in
\mathbb Z}\left(\mathbf
E|\eta_{jk}|^2\right)^{1/2}\cdot2^{j/2}C_{\psi}\right)^p,\label{cn+}\end{eqnarray}
where $J:=\min\{n,\min\{j\in \mathbb{N}_0:\, k_j<2^jT+1\}\}.$
\end{thm}
\noindent{\it Proof.} Notice that $J\to \infty,$ when $\mathbf k_n\to\infty.$ By the choice of $J$ and (\ref{supsumksieta^2}) we get $c_n^\infty\to 0,$ when $\mathbf k_n\to\infty.$ The statement of the theorem  follows from
Theorem~\ref{conv_Lp_Sub_ner}, Corollary~\ref{cor2_str_Sub},
Lemma~\ref{lem_ner_phi}, and the fact $\sum_{k\in \mathbb
Z}\psi_{jk}(t)=2^{j/2}\sum_{k\in \mathbb Z}\psi(2^jt-k).$ \hfill   $\Box$ 

\begin{thm}\label{str_Sub_covLp4}
Let $\mathbf{X}(t),$ $t\in\mathbb R,$ be a strictly
$Sub_{\varphi}(\Omega)$ stationary  short-memory random process  with the covariance function $R(t-s):=\mathbf
E\mathbf{X}(t)\overline{\mathbf{X}(s)}.$ Let the assumptions of
 Theorem~{\rm\ref{8103}} are satisfied and
$\int_0^T\left(\mathbf E(\mathbf X(t))^2\right)^{p/2}\,dt<+\infty.$ Suppose that there exists  the Fourier transform $\widehat{R}(z)=\int_{\mathbb R}e^{-izt}{R(t)}dt$ and for some
$\alpha>0:$
$$\int_{\mathbb
R}\left|\widehat{R}(z)\right|\cdot\left|z\right|^{\alpha}\,dz<\infty.$$
If $\widehat{\psi}(\cdot)$ is a Lipschitz function of order
$\alpha/2,$  then {\rm(\ref{supsumksieta^2})} holds.
\end{thm}
\noindent{\it Proof.} By Parseval's theorem and the representation $\widehat{\psi}_{jk}(z)=2^{-\frac
j2}e^{-i \frac k{2^j} z}\widehat{\psi}\left(z/{2^j}\right),$ we deduce
$$\mathbf E|\eta_{jk}|^2=\left|\int_{\mathbb R}\int_{\mathbb R}R(u-v)\overline{\psi_{jk}(u)}\,du\,\psi_{jk}(v)\,dv\right|=\left|\int_{\mathbb R}\int_{\mathbb R}\frac{e^{-ivz}}{2 \pi}\widehat{R}(z)\overline{\widehat{\psi}_{jk}(z)}\,dz\psi_{jk}(v)\,dv\right|$$
\begin{equation}\label{st}\le\frac{1}{2 \pi}\int_{\mathbb R}\left|\widehat{R}(z)\right|\,\left|
\overline{\widehat{\psi}_{jk}(z)}\,\widehat{\psi}_{jk}(z)\right|\,dz=\frac{1}{2^{j+1}
\pi}\int_{\mathbb
R}\left|\widehat{R}(z)\right|\cdot\left|\widehat{\psi}\left(\frac
z{2^j}\right)\right|^2\,dz.
\end{equation}
By the Lipschitz conditions, (\ref{st}), and $\widehat{\psi}(0)=0,$ we obtain that  for all
$k\in\mathbb Z:$
\begin{equation}\label{st1}\left|\mathbf
E\eta_{jk}\right|^2\le\frac{C}{2^{j+1} \pi}\int_{\mathbb
R}\left|\widehat{R}(z)\right|\cdot\left(\frac
{|z|}{2^j}\right)^{\alpha}\,dz\le\frac{C}{\pi\,
2^{1+j(1+\alpha)}}\int_{\mathbb
R}\left|\widehat{R}(z)\right|\cdot\left|z\right|^{\alpha}\,dz,\end{equation}
where $C>0$ is the Lipschitz constant.

Similarly to (\ref{st}) we deduce
\begin{equation}\label{st2}\mathbf E|\xi_{0k}|^2\le \frac{1}{2\pi}\int_{\mathbb
R}\left|\widehat{R}(z)\right|\cdot\left|\widehat{\phi}\left(z\right)\right|^2\,dz<\infty\,.
\end{equation}
The integral in (\ref{st2}) is finite because the Assumption S implies $\phi(\cdot)\in L_1(\mathbb R)$ and therefore $\widehat{\phi}(\cdot)$ is bounded.

By  estimates (\ref{st1}) and (\ref{st2}) the series in (\ref{supsumksieta^2}) is convergent.
\hfill   $\Box$
\begin{thm}\label{str_Sub_covLp_nonstation}
Let $\mathbf{X}(t),$ $t\in\mathbb R,$ be a strictly
$Sub_{\varphi}(\Omega)$ random process
 with the continuous
covariance function $R(t,s):=\mathbf
E\mathbf{X}(t)\overline{\mathbf{X}(s)}.$  Suppose that
$\int_0^T\left(\mathbf E(\mathbf X(t))^2\right)^{p/2}\,dt<+\infty,$  there exist
$\widehat{R}_2(z,w):=\int_{\mathbb R}\int_{\mathbb
R}e^{-izu}e^{-iwu}{R(u,v)}\,dudv,$ and for some $\alpha>0:$ 
$$\int_{\mathbb
R}\int_{\mathbb
R}\left|\widehat{R}_2(z,w)\right|\cdot\left|z\right|^{\alpha}\cdot\left|w\right|^{\alpha}\,dzdw<\infty.$$
If  the assumptions of
 Theorem~{\rm\ref{8103}} are satisfied and $\widehat{\psi}(\cdot)$ is a Lipschitz function of order~$\alpha,$ then  {\rm (\ref{supsumksieta^2})} holds true.
\end{thm}
\noindent{\it Proof.} By Parseval's theorem,
\begin{eqnarray}
\mathbf
E|\eta_{jk}|^2&=&\int_{\mathbb R}\int_{\mathbb R}\mathbf E
\mathbf{X}(u)\overline{\mathbf{X}(v)}\
\overline{\psi_{jk}(u)}\psi_{jk}(v)\,dudv\nonumber\\
&=&\frac{1}{(2
\pi)^2}\int_{\mathbb R}\int_{\mathbb
R}\widehat{R}_2(z,w)\,\overline{\widehat{\psi}_{jk}(z)}\,\widehat{\psi}_{jk}(w)\,dz\,dw\nonumber\\
&\le&\frac{1}{2^{j}(2 \pi)^2}\int_{\mathbb R}\int_{\mathbb
R}\left|\widehat{R}_2(z,w)\right|\cdot\left|\widehat{\psi}\left(\frac
z{2^j}\right)\right|\cdot\left|\widehat{\psi}\left(\frac
w{2^j}\right)\right|\,dz\,dw.\label{st1_2}
\end{eqnarray}

 By properties of the $m$-wavelet
$\psi(\cdot)$ we have $\widehat{\psi}(0)=0.$ Therefore, using the Lipschitz
conditions, we obtain
$$\mathbf
E|\eta_{jk}|^2\le\frac{C^2}{(2 \pi)^2 2^{j(1+2\alpha)}}\int_{\mathbb
R}\int_{\mathbb
R}\left|\widehat{R}_2(z,w)\right|\cdot\left|z\right|^{\alpha}\cdot\left|w\right|^{\alpha}\,dz\,dw,$$
where $C>0$ is the Lipschitz constant.

Analogously to (\ref{st1_2})  we deduce that
\[\mathbf E|\xi_{0k}|^2\le
\frac{1}{4 \pi^2}\int_{\mathbb R}\int_{\mathbb
R}\left|\widehat{R}_2(z,w)\right|\cdot\left|\widehat{\phi}\left(
z\right)\right|\cdot\left|\widehat{\phi}\left(
w\right)\right|\,dz\,dw<\infty
\]
due to the boundedness of $\widehat{\phi}(\cdot).$

An  application of the above estimates gives the convergence of (\ref{supsumksieta^2}).
 \hfill   $\Box$
 
 \begin{rmk}  Conditions of Theorems {\rm\ref{str_Sub_covLp4}} and
  {\rm\ref{str_Sub_covLp_nonstation}} on the random process $\mathbf{X}(t)$ are formulated in terms of its spectral density. These conditions are related to the behavior of the high-frequency part of the spectrum. Such assumptions are standard in the convergence studies of stochastic approximations.
 \end{rmk}
\begin{rmk} If the assumptions of Theorem {\rm\ref{str_Sub_covLp4}} or
{\rm\ref{str_Sub_covLp_nonstation}} are satisfied, then
 {\rm(\ref{P_int_Xn-Xm<exp})} holds true for $c_n^\infty$ given by {\rm(\ref{cn+})}. All the terms in {\rm(\ref{cn+})} can be easily computed in practice for specific stochastic processes and wavelet bases.
\end{rmk}

\section{Conclusions}
The obtained results may have various practical applications for the approximation and simulation of random processes. The analysis of the rate of convergence provides a constructive algorithm for determining the number of terms in the wavelet expansions to ensure the approximation of stochastic processes with given accuracy.

The developed methodology and results are important extensions of the recent findings in the wavelet approximation theory of stochastic processes to the space $L_p([0,T])$ and the class of $\varphi$-sub-Gaussian random processes. This class plays an important role in generalizations of various theoretical properties of Gaussian processes.
In addition to classical applications of $\varphi$-sub-Gaussian random processes in signal processing, the results can also be used in new areas, like compressed sensing and actuarial modelling, consult, for example,  \citealt{lab,ver,yam}.

It would be of interest 
\begin{itemize}
\item to adopt and specify the results for different wavelet bases, which satisfy the assumptions of the theorems, for example, Daubechies, Battle-Lemarie and Meyer wavelets;
\item to examine the tightness of the estimates by simulations.
\end{itemize}


\begin{thebibliography}{}

\bibitem[Atto, Berthoumieu(2010)]{att} Atto  A, Berthoumieu Y (2012)
Wavelet packets of nonstationary random processes: contributing factors for stationarity and decorrelation. 
IEEE Trans Inform Theory 58(1):317–330

\bibitem[Bardet, Tudor(2010)]{bar} Bardet  JM, Tudor CA (2010)
A wavelet analysis of the Rosenblatt process: chaos expansion and estimation of the self-similarity parameter. 
Stochastic Process Appl 120(12):2331–2362

\bibitem[Buldygin, Kozachenko(2000)]{bul} Buldygin VV,  Kozachenko YuV (2000) Metric
characterization of random variables and random processes.
 American Mathematical Society, Providence R.I.

\bibitem[Cambanis, Masry(1994)]{cam} Cambanis S,   Masry E (1994) Wavelet
approximation of deterministic and random signals: convergence
properties and rates. IEEE Trans Inf Theory\linebreak
40(4):1013-1029

\bibitem[Clausel et al.(2012)]{cla} Clausel M, Roueff F, Taqqu MS, Tudor C (2012) Large scale behavior of wavelet coefficients of non-linear subordinated processes with long memory. Appl Comput Harmon Anal 32(2): 223–241

\bibitem[Daubechies(1992)]{dau} Daubechies  I (1992) Ten lectures on
wavelets. SIAM, Philadelphia

\bibitem[ Didier, Pipiras(2008)]{did} Didier G,   Pipiras V (2008) Gaussian stationary processes: adaptive wavelet decompositions, discrete approximations and their convergence. J Fourier Anal and Appl 14:203-234

\bibitem[Giuliano Antonini et al.(2003)]{ant} Giuliano Antonini R, Kozachenko YuV,
Nikitina T (2003) Spaces of $\varphi$-sub-Gaussian random
variables. Mem Mat Appl 121(27)\,fasc~1:95-124

\bibitem[Hardle et al.(1998)]{har} Hardle  W, Kerkyacharian G, Picard  D,
Tsybakov A (1998)  Wavelets, approximation and
statistical applications.  Springer, New York

\bibitem[Jaffard(2001)]{jaf} Jaffard S (2001) Wavelet expansions,
function spaces and multifractal analysis. In:  Byrnes JS (ed)
Twentieth century harmonic analysis -- a celebration.
 Kluwer Acad Publ, Dordrecht, pp 127-144
 
\bibitem[Kurbanmuradov, Sabelfeld(2008)]{kur} Kurbanmuradov O, Sabelfeld K  (2008) Convergence of fourier-wavelet models for Gaussian random processes. SIAM J Numer Anal 46(6):3084-3112
 
 \bibitem[Kozachenko, Kamenshchikova(2009)]{koz0}  Kozachenko Yu, Kamenshchikova O (2009) Approximation of $\operatorname {SSub}_{\varphi }(\Omega )$ stochastic processes in the space $L_{p}(\mathbb {T})$. Theor Probability and Math Statist 79:83-88

\bibitem[Kozachenko, Kovalchuk(1985)]{koz1} Kozachenko Yu, Kovalchuk Yu (1998)
Boundary value problems with random initial conditions and
functional series from $Sub_{\phi}( \Omega)$. Ukrainian
Math J 50:504-515

\bibitem[Kozachenko et al.(2011)]{koz2} Kozachenko Yu,  Olenko A,
Polosmak O (2011) Uniform convergence of wavelet expansions of
Gaussian random processes. Stoch Anal Appl 29:169-184

\bibitem[Kozachenko et al.(2013)]{koz3} Kozachenko Yu,  Olenko A,
Polosmak O (2013) Convergence rate of wavelet expansions of Gaussian random processes. will appear in Comm Statist Theory Methods 

\bibitem[Kozachenko, Ostrovskyi(1985)]{koz4} Kozachenko Yu, Ostrovskyi E (1985) Banach
spaces of random variables of Sub-gaussian type. Theor
Probability and Math Statist 32:42-53

\bibitem[Kozachenko, Polosmak(2008)]{koz5} Kozachenko Yu, Polosmak O (2008)
 Uniform convergence in probability of wavelet expansions of
random processes from $L_2( \Omega)$. Random Oper Stoch Equ\linebreak 16(4):12-37

\bibitem[Labate et al.(2013)]{lab}  Labate  D,  Weiss G,   Wilson E (2013) Wavelets. Notices Amer Math Soc 60(1):66-76 

\bibitem[Triebel(2008)]{tri}  Triebel H (2008) Function spaces and
wavelets on domains.  European Mathematical Society, Zu\"{u}rich

\bibitem[Vershynin(2012)]{ver}   Vershynin R (2012) Introduction to the non-asymptotic analysis of random matrices. In  Eldar Y, Kutyniok G (ed)  Compressed sensing, theory and applications. Cambridge University Press, Cambridge, pp 210-268

\bibitem[Yamnenko (2006)]{yam}   Yamnenko R (2006) Ruin probability for generalized $\varphi$-sub-Gaussian fractional Brownian motion. Theory Stoch Process 12(28)1-2:261-275

\bibitem[Zhang, Waiter(1994)]{zha} Zhang J,  Waiter G (1994) A wavelet-based KL-like expansion for wide-sense stationary random processes. IEEE Trans Signal Proc 42(7):1737-1745
\end{thebibliography}
\end{document}